\documentclass[lettersize,journal]{IEEEtran}
\usepackage{amsmath,amsfonts}
\usepackage{algorithmic}
\usepackage{algorithm}
\usepackage{array}
\usepackage[caption=false,font=normalsize,labelfont=sf,textfont=sf]{subfig}
\usepackage{textcomp}
\usepackage{stfloats}
\usepackage{url}
\usepackage{verbatim}
\usepackage{graphicx}
\usepackage{cite}
\usepackage{hyperref}
\usepackage{array}
\usepackage{titlesec}
\titlespacing*{\subsection}{0pt}{1ex}{1ex}
\titlespacing*{\section}{0pt}{1ex}{1ex}
\newtheorem{assumption}{Assumption}

\begin{document}

\title{Toward Robust Manufacturing Scheduling: Stochastic Job-Shop Scheduling}

\author{Mikhail A. Bragin,~\IEEEmembership{Senior Member,~IEEE}, Matthew E. Wilhelm, Nanpeng Yu,~\IEEEmembership{Senior Member,~IEEE} and Matthew D. Stuber
\thanks{Mikhail A. Bragin is with 
the Institute for Advanced Systems Engineering, University of Connecticut, Storrs, CT, 06269, USA; Matthew E. Wilhelm and Matthew D. Stuber are with the Process Systems and Operations Research Laboratory, Department of Chemical and Biomolecular Engineering, University of Connecticut, Storrs, CT, 06269, USA; Nanpeng Yu is with the Department of Electrical and Computer Engineering, University of California Riverside, Riverside, CA, 92521, USA;}
}

\markboth{Journal of \LaTeX\ Class Files,~Vol.~14, No.~8, August~2021}%
{Shell \MakeLowercase{\textit{et al.}}: A Sample Article Using IEEEtran.cls for IEEE Journals}


\maketitle

\begin{abstract}
Manufacturing plays a significant role in economic development, production, exports, and job creation, which ultimately contribute to improving the quality of life. The presence of manufacturing defects is, however, inevitable leading to products being discarded, i.e. scrapped. In some cases, defective products can be repaired through rework. Scrap and rework cause a longer completion time, which can contribute to orders being shipped late. Moreover, the presence of uncertainties and combinatorial complexity significantly increases the difficulty of complex manufacturing scheduling. This paper tackles this challenge, exemplified by a case study on stochastic job-shop scheduling in low-volume, high-variety manufacturing contexts. To ensure on-time delivery, high-quality solutions are required, and near-optimal solutions must be obtained within strict time constraints to ensure smooth operations on the job-shop floor. To efficiently solve the stochastic job-shop scheduling (JSS) problem, a recently-developed Surrogate ``Level-Based'' Lagrangian Relaxation is used to reduce computational effort while efficiently exploiting the geometric convergence potential inherent to Polyak's step-sizing formula thereby leading to fast convergence. Numerical testing demonstrates that the new method is two orders of magnitude faster as compared to commercial solvers. 

\textit{Note to Practitioners---}Manufacturing defects leading to scrap or rework create significant challenges, from both computational as well as on-time delivery standpoints. To assist practitioners in overcoming these challenges, this paper presents a case study of stochastic job-shop scheduling as well as Surrogate ``Level-Based" Lagrangian Relaxation as an efficient solution methodology, which not only reduces the time required to obtain high-quality, near-optimal solutions for complex scheduling problems but also minimizes the dependence on hyperparameters. As a result, the method becomes more user-friendly, reducing the need for domain knowledge or instance-specific heuristical adjustments. By implementing the proposed method, practitioners in the manufacturing industry can expect a substantial increase in the speed of obtaining solutions to their scheduling problems, expecting to significantly outperform commercial solvers. This improvement will help organizations meet delivery deadlines, maintain high-quality production standards, and ultimately enhance their competitiveness in the market.
\end{abstract}

\begin{IEEEkeywords}
Advanced Manufacturing, Stochastic Job-Shop Scheduling, Lagrangian Relaxation, Markov Processes
\end{IEEEkeywords}

\section{Introduction}
\IEEEPARstart{M}{manufacturers'} ultimate goal is customer satisfaction, achieved through delivering high-quality products on time to remain competitive in the marketplace \cite{cheng15,choi15,rahman21}. 
Late deliveries, however, can diminish the value of products in the eyes of customers. Agile and flexible tailored production greatly influence the value created. The trend towards improving flexibility is through make-to-order production, allowing the ordering of bespoke products \cite{chen15, rahman17, wang17}. Manufacturing operations must account for tardiness to avoid late shipments and subsequent tarnishing reputations in the marketplace. Correcting quality problems at the factory level is more cost-effective than resolving quality problems after the product has reached the end customer \cite{ergun2022structured}. To this end, machining-level actions (e.g., feed rate and depth of cut) \cite{wilson2022multi}, as well as tardiness and reactive re-scheduling, have been optimized \cite{sun2022simulation}.

\textbf{Scope:} This paper addresses tardiness-based job-shop scheduling (JSS), which is crucial 
to avoid late shipments and maintain customer satisfaction. Job shops are a specific type of manufacturing environment designed for low-volume/high-variety production, making them ideal for producing bespoke products. However, the JSS problem is 
difficult due to the combinatorial complexity caused by discrete scheduling decision variables.
As a result, finding even near-optimal solutions can be challenging, leading to inefficient operations and excessive costs, or long solving times that contribute to late deliveries and associated consequences.

Additionally, the paper considers the possibility of damaged parts or orders resulting in scrapped or reworked parts. The problem is modeled stochastically, leading to a probabilistic routing of jobs that takes uncertainties into account. 
Inaccurate decisions in the presence of uncertainty can be particularly costly, as scrapping one operation of one job may require remanufacturing the entire part, causing unexpected delays in production that can lead to customer dissatisfaction.
This paper aims to obtain near-optimal solutions quickly to mitigate the impact of these uncertainties and delays.

The rest of this section is organized as follows: In Subsection \ref{sec:IntMinTardiness}, we review JSS problems with the consideration of tardiness. In Subsection \ref{sec:IntStochJSS}, we review stochastic JSS. In Subsection \ref{sec:IntJSSform}, we review recent formulations of JSS problems. Finally, in Subsection \ref{sec:IntJSSsol}, we provide a review of solution methodologies for solving JSS problems and general discrete programming problems.
\subsection{Reducing Tardiness: On-Time Delivery in Job Shops}
\label{sec:IntMinTardiness}
The goal of a JSS problem is to assign a set of jobs, each consisting of a sequence of consecutive \textit{non-preemptive} (uninterruptible)
operations, to be processed on a set of machines. Generally, JSS problems have been recognized as difficult combinatorial optimization problems \cite{garey76, zhu20}.

As previously reviewed, on-time delivery has long been recognized as an important aspect of production; the minimization of tardiness has been used as an objective within JSS problems as well. The seminal works on ``tardy problems'' begin with the research of \cite{sidney77} and \cite{lakshminarayan78} with subsequent research of, for example, \cite{vepsalainen87}, \cite{chryssolouris88}, \cite{anderson90}, \cite{du90}, \cite{chryssolouris91}, \cite{kanet93}, and \cite{rohleder93}. 
The problem is known to be NP-hard (e.g., \cite{du90}), even single-machine scheduling \cite{seidmann81, panwalkar82, ow89,du90} belongs to a class of NP-hard problems. Yet, the job-shop scheduling problems are also subject to strict computational requirements; long-solving times are not allowed since they may contribute to orders shipped late.

For this reason, several rules have been developed in early implementations such as ``first-come, first-serve'', ``shortest processing time'', or ``earliest operation due date'' rules; these and several other rules are described and compared by \cite{kanet93}. The minimization of tardiness, however, requires the coordination of multiple jobs and multiple machines, while dispatching rules typically ignore other resources.  In fact, as pointed out by \cite{barman98}, ``no single priority rule performs well.'' Therefore, most of the research has been heavily invested in the development of solution methodologies beyond rule-based ones.
\subsection{Robustness to Uncertainty: Stochastic Job-Shop Scheduling}
\label{sec:IntStochJSS}
Uncertainty is a concomitant process within any production system and must be accounted for. 
Therefore, it is important to account for potential causes of disruptions to mitigate the consequences should such unforeseen events occur in actual operations. Without the appropriate incorporation of uncertainty into problem formulations, deterministic scheduling decisions may not be optimal during operations, if feasible at all. In terms of JSS, 
neglecting uncertainties within the optimization may result in unexpectedly high tardiness in actual operations, thereby compromising the on-time delivery of products.

In the extant research, uncertain processing times have been considered \cite{luh99, Golenko02, Lei07, lei11, Horng12, zhang13, yang14, shen16, jamili19, horng21}\footnote{Exceptions are, perhaps, \cite{Hasan11, Lin19} where machine unavailability is considered, though, both papers consider minimization of makespan, not tardiness.}. Within these papers, while the processing time is uncertain, a job is assumed to be processed with certainty. However, in actual operations, a product may be damaged (scrapped) after being processed and must be remade anew thereby potentially leading to higher tardiness. This addition of tardiness cannot be captured through stochastic modeling of processing times, since a job may require several operations, all of which may need to be repeated. Another salient feature of the above existing methods is their meta-heuristic nature (with the exception of the work by \cite{luh99}). Heuristic-based algorithms generally do not provide a lower bound to quantify the solution quality. Some notable exceptions are algorithms such as those proposed by \cite{dawande06}.

A fast resolution of JSS problems rests upon two pillars: formulation and solution methodology, both are reviewed next. 
\subsection{``Tight'' Job-Shop Scheduling ILP/MILP Formulations}
\label{sec:IntJSSform}
An aspect of problem formulation that has a direct impact on the solution process is the formulation ``tightening,'' the gist of which is to delineate facets of the convex hull. If successful, then the combinatorial problem reduces to a much easier-to-solve LP problem; even if a problem is ``partially'' tightened, the CPU time is greatly reduced \cite{yan18, yan21}. A recent advancement reduces decision variables and constraints in ILP formulations, proving to be tighter than previous ones \cite{Liu21}: only the beginning times of operations are decision variables thereby reducing the number of decision variables and constraints compared to those within the previous formulations. Since the formulation has few decision variables and constraints, it leads to several orders of magnitude improvement over previously-used approaches in terms of CPU time and to optimality of solution obtained (from 3,600 s and 3.72\% gap down to 3.31 s and 0\% gap) for a problem instance from \cite{Hoitomt93}, when solved by using standard B\&C \cite{Liu21}.
\subsection{Previous Solution Methodologies}
\label{sec:IntJSSsol}
In this subsection, methods used to solve general scheduling problems are briefly reviewed. Then, methods specifically developed for the JSS problem are reviewed. Since the methodology of this paper is based on Lagrangian Relaxation, the LR-based methods are reviewed as related to JSS problems. Finally, the most recent advancements of LR are presented. 

\noindent \textbf{Methodologies for General Scheduling Problems.} To solve scheduling problems (including JSS problems minimizing makespan, and flow shop scheduling), the following methods have been used: ant colony optimization \cite{leung10}, conflict-directed search \cite{grimes2015solving}, constraint programming \cite{beck2011combining, malapert2012optimal}, decomposition methods \cite{zhao2019decomposition}, various version of evolutionary algorithms \cite{kim03, wu2019moels}, greedy algorithms \cite{zhao2021iterated, zhao2020heuristic}, various versions of genetic algorithms \cite{li16, li18, zhu20, li2021scoring, cao2022two}, heuristics/metaheuristics \cite{dawande06, zhao2020heuristic, pan2022improved}, Lagrangian Relaxation \cite{ilo18, 10097555}, local search \cite{aarts1994computational, vaessens1996job, beck2011combining}, particle swarm optimization \cite{nouiri18}, Petri nets \cite{casalino2019optimal, zhao2020heuristic}, reinforcement learning \cite{park2019reinforcement, luo2021real, song2023stochastic, li2023scheduling} and taboo search \cite{taillard1994parallel}.    
 

\noindent \textbf{Methodologies for Job-Shop Scheduling Problems with the Consideration of Tardiness.} Because of the complexity, JSS problems require a more sophisticated methodology beyond the simple rules reviewed in subsection \ref{sec:IntMinTardiness}. Toward this goal, several methods have been developed such as Lagrangian Relaxation \cite{Hoitomt93, luh99}, genetic algorithms \cite{herrmann95, lei11}, heuristics \cite{Golenko02}, neural networks with simulated annealing \cite{Tavakkoli05}, evolutionary algorithms \cite{Lei07, Horng12, zhang13}, ordinal optimization \cite{yang14}, and evolutionary learning-based simulation optimization \cite{Ghasemi21}. 

\noindent \textbf{Lagrangian Relaxation for Job-Shop Scheduling Problems.} To address the combinatorial complexity inherent in Job Shop Scheduling (JSS) problems, Lagrangian Relaxation (LR) has been employed to ``reverse" the complexity, leading to an exponential reduction in the effort required to solve subproblems. For JSS problems, relaxing machine capacity constraints allows the subproblems associated with individual jobs to be solved with significantly reduced complexity. The application of standard Lagrangian Relaxation in JSS problems can be traced back to the work of Hoitomt et al. \cite{Hoitomt93}, where Lagrange multipliers were updated using the Polyak stepsize \cite{Polyak69} and sub-gradient directions. While there were other advancements in improving the coordination of LR, in particular, to solve JSS problems \cite{Bragin15, yan18, Bragin19, Liu21}, Polyak's results have the potential to achieve geometric convergence (fastest possible) and deserve a separate discussion. 



\noindent \textbf{Surrogate Sub-Gradient Method.} Within the Surrogate Sub-Gradient Method \cite{zhao99}, multipliers are updated after solving one subproblem at a time rather than solving all the subproblems as within the standard sub-gradient method. This significantly reduces the computational effort. The step-sizes are updated by using the following variation of the Polyak's formula (originally developed in \cite{Polyak69}): 
\vspace{-2mm}
\begin{align}
& 0 < s^k < \gamma \cdot \frac{q(\lambda^{*}) - L^k}{\left\|\tilde{g}^k\right\|^2_2}, \gamma < 1, \label{SSM}
\end{align}
where $L^k$ is the Lagrangian value and $\tilde{g}^k$ are the levels of constraint violations (``surrogate'' sub-gradients). These values are used in place of dual values $q(\lambda^k)$ and sub-gradients $g^k$.  The convergence to $\lambda^{*}$ is guaranteed \cite{zhao99}\footnote{Unlike that in Polyak's formula \cite{Polyak69}, parameter $\gamma$ is less than 1 rather than 2.}. The concomitant reduction of multiplier zigzagging has been also observed. Even though the \textit{geometric} convergence rate is expected, e.g., multipliers strictly move closer to $\lambda^*$, the implementation of the method is problematic because of the unavailability of the knowledge about the optimal dual value $q(\lambda^*)$. 

\noindent \textbf{Surrogate ``Level-Based'' Lagrangian Relaxation (SLBLR).} To overcome the above issue, the SLBLR method has been recently developed by \cite{bragin2022surrogate} to obtain “level” (over)estimates of the optimal dual value $q(\lambda^*)$ within the Polyak's step-sizing formula \eqref{SSM} to exploit the geometric convergence potential. 

The main idea of SLBLR is: if multipliers do not approach $\lambda^*$, then stepsizes violate \eqref{SSM} as proved in \cite{bragin2022surrogate}. Specifically, if the following ``multiplier-divergence-detection'' problem (with $\lambda$ being a continuous decision variable: $\lambda \in \mathbb{R}^M$)
\begin{equation} 
    \begin{cases}
      \|\lambda-\lambda^{k_j+1}\|_2 \leq \|\lambda-\lambda^{k_j}\|_2,\\
      \|\lambda-\lambda^{k_j+2}\|_2 \leq 
\|\lambda-\lambda^{k_j+1}\|_2,\\
\qquad\qquad\vdots\\
\|\lambda-\lambda^{k_j+n_j}\|_2 \leq 
\|\lambda-\lambda^{k_j+n_j-1}\|_2, \label{lambdfeas}
    \end{cases}
\end{equation}
admits no feasible solution with respect to $\lambda$ for some $k_j$ and $n_j$, then there exists $i \in [k_j,k_j+n_j]$ such that 
\begin{equation}\label{violstep}
s^i \ge \gamma \cdot \frac{q(\lambda^{*}) - L^k}{\|\tilde{g}^k\|^2_2}. 
\end{equation} 

From \eqref{violstep} it follows that there exists an overestimate of $q(\lambda^{*})$, and the overestimate (the ``level value'') has been derived to be: 
\begin{equation} 
\overline{q}_{j} = \max_{i \in [k_j,k_j+n_j]} \left(\frac{1}{\gamma} \cdot s^i \cdot \|\tilde{g}^k\|^2_2 + L^k\right) > q(\lambda^{*}). \label{eq28}
\end{equation}

As a result, the step-sizes are set by using a series of decreasing (not necessarily monotonically) overestimates $\{\overline{q}_{k}\}$ as:  
\begin{equation}
 s^k = \zeta \cdot \gamma \cdot \frac{\overline{q}_{k} - L^k}{\left\|\tilde{g}^k\right\|^2_2}, \gamma < 1, \zeta < 1.  \label{SLBLR}
\end{equation}
The factor of $\gamma$ in the above formula is inherited from the formula \eqref{SSM}, and the extra factor of $\zeta$ is introduced to counter the fact that the stepsizes are set by using overestimates of the optimal dual value $q(\lambda^*)$. 

The rest of the paper is organized as follows. After presenting stochastic JSS formulation in Section \ref{sec:JSSform}, the SLBLR method is used to solve the problem in Section \ref{sec:SolMethod}. 
The numerical case studies and solution results are presented in 
Section \ref{sec:numerical_examples} demonstrates the advantages of the SLBLR method. The conclusion is provided in Section \ref{sec:conclusion}.

\section{Job-Shop Scheduling Problem Description and Formulation}
\label{sec:JSSform}

Within a typical job shop, machines are grouped by their functionality (e.g., milling, lathing, drilling, etc.), and the number of such groups is denoted by $M$. Within each group, the number of machines referred to as the ``$\underline{\text{c}}$apacity'' is denoted as $C_m,$ where $m = 1,\dots,M.$ For every order $i = 1,\dots, I$ (e.g., a part to be processed) arriving at a job shop with priority $w_i$
a specific sequence of $J_i$ operations, processing times $p_{i,j,m}$ (for each part $i$, operation $j$\footnote{Operation $j$ is denoted as $(i,j)$ whenever appropriate to indicate specific part $i$ requiring operation $j$.} and machine type $m$), as well as an integer-valued due date $d_i$ is specified. Processing times are assumed to be an integer and each operation can occupy several contiguous time blocks, that is, every operation is assumed to be \textit{non-preemptive}. Accordingly, the time horizon is assumed to consist of $T$ discrete equidistant time blocks. A part $i$ may need to go through a sequence of operations, each is eligible to be performed by one or a few machine groups. Let $E_{i,j}$ be a group of machines $\underline{\text{e}}$ligible to process operation $(i,j)$, and $O_m$ be a set of $\underline{\text{o}}$perations that can be processed on machine $m$. 

\begin{figure}[ht]
\centering
\includegraphics[height=2.25in]{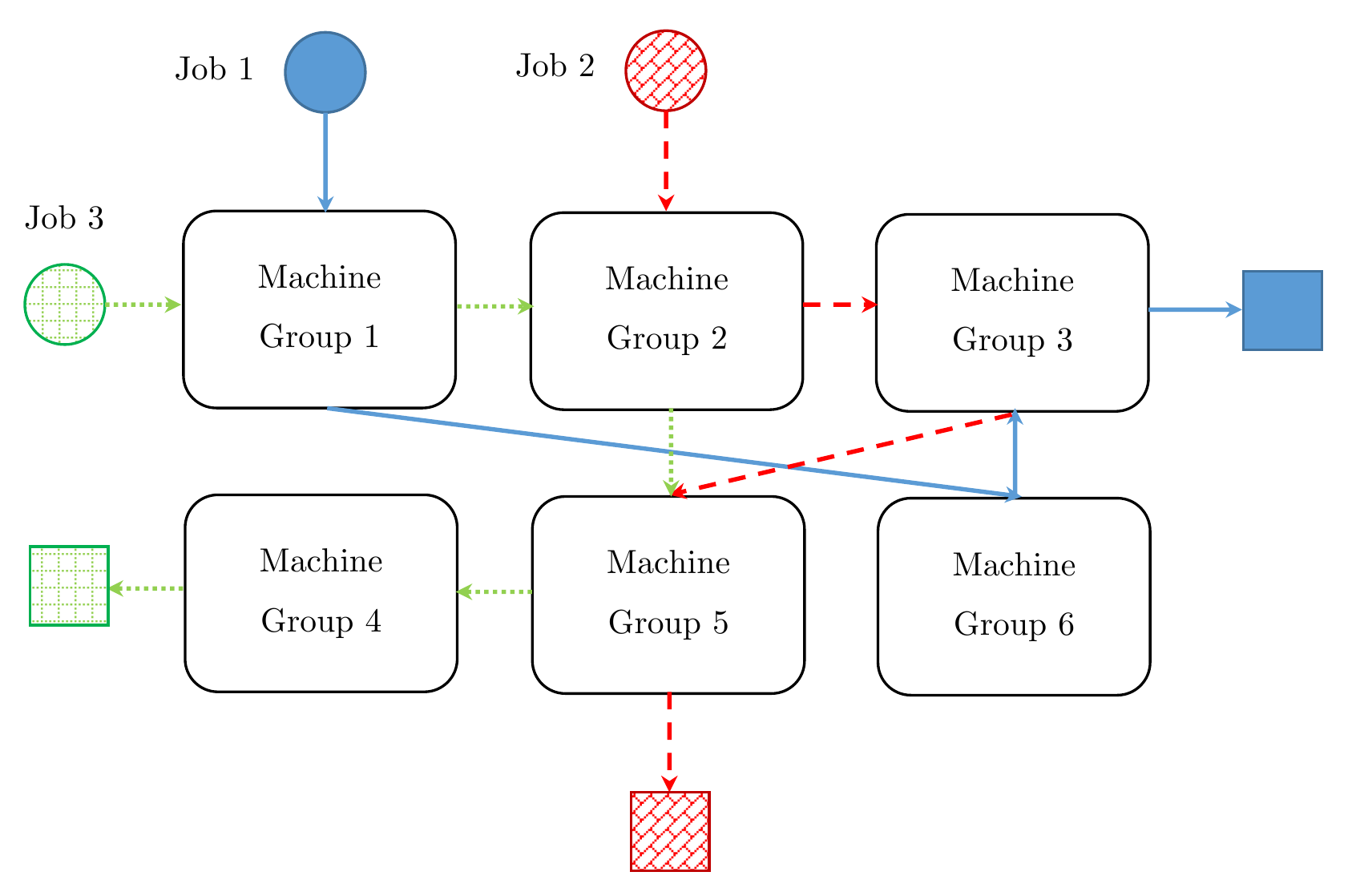}
\caption{A schematic of job flow within a typical job shop. Circles represent inputs for each job and squares represent the corresponding products as outputs. Machine groups may contain 1 or more machines.} 
\label{frontier}
\end{figure}

Figure \ref{frontier} illustrates possible, albeit simplified, job flow within a job shop omitting the time aspect. 
As shown in Figure \ref{frontier}, Job 1 requires 3 operations to be performed by Machine Groups 1, 6, and 3, sequentially. Likewise, Job 2 requires 3 operations to be performed by Machine Groups 2, 3, and 5, sequentially; and Job 3 requires 4 operations to be performed by Machine Groups 1, 2, 5, and 4, sequentially. If Machine Group 1 contains only one machine available and Jobs 1 and 3 arrive at the same time, only one job can be processed during the next time period --- this is to be decided through optimization. A similar argument holds for other machine groups, except Machine Groups 4 and 6, which only 
operates on Job 3 and Job 1, respectively. In terms of the job shop shown in Figure \ref{frontier}, $E_{1,1} = \{1\}$, $E_{1,2} = \{6\}$, etc., and $O_1 = \{(1,1),(3,1)\}, O_2 = \{(2,1),(3,2)\},$ etc. 

\begin{figure}[ht]
		\centering
 \includegraphics[trim = 0 0 0 0, scale=0.35, angle=0]{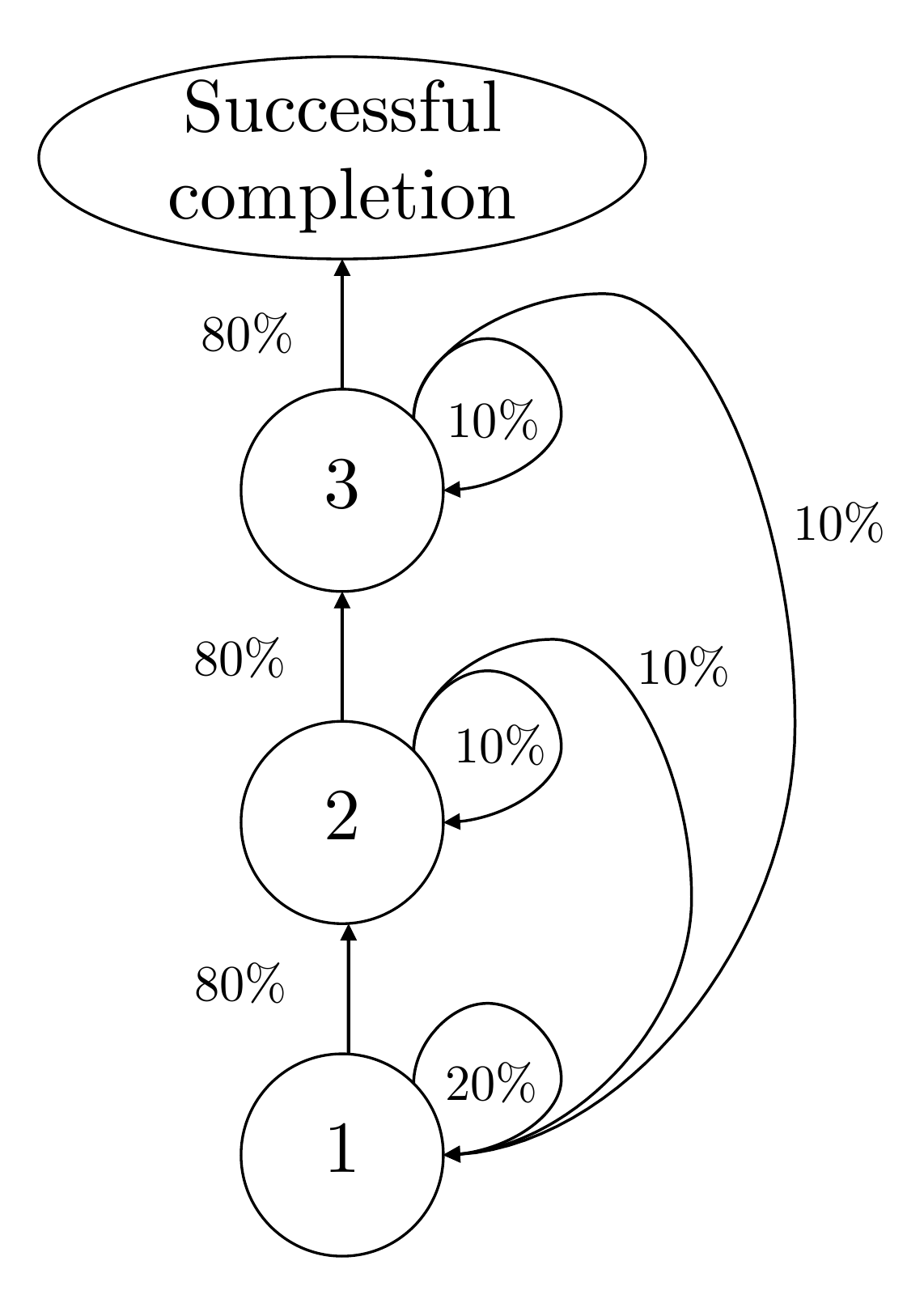}		
 \caption{A Markov process is illustrated and used in this work to capture transitions among states (operations) of a part. In this illustration, with a probability of 10\% a part is either scrapped or reworked; at Operation 1, scrap and rework are lumped together since a defective part needs to begin getting processed at Operation 1.}
 \label{fig:markov}
\end{figure}

There is a nonzero probability $p_{(i,j)}^s$ of scrap after operation $j$ of part $i$, and production must restart from Operation $1$. Additionally, there is a nonzero probability $p_{(i,j)}^r$ that a scrapped part $i$ can be recovered through rework and its production must restart from operation $j$. Transitions among the states of a part follow a Markov process schematically illustrated in Figure \ref{fig:markov}. For simplicity of exposition, probabilities of scrap and rework are chosen to be 20\% and 50\%, respectively. Transitions from one operation to the next thus happen with a probability of 80\%, and in the case of a defective part, with equal probabilities of 10\%, a part will either be reworked or completely discarded, triggering the restart of processing from Operation 1.

\begin{assumption}\label{assume1}
The models are developed in this paper under the following assumptions:

\begin{enumerate} 
    \item Each machine can process at most one operation at a time;\label{a1}
    \item Preemption is not allowed, so no operating can be interrupted during processing;
    \item No two operations of a given job can be processed at the same time, and moreover, the predefined sequence of the operations should be kept: no two operations can be swapped;\label{a3}
    \item All jobs and machines are available at Time 1;
    \item All the non-processing related times such as setup times and transportation times are included in the processing time.
\end{enumerate}
\end{assumption}

The job-shop scheduler's task is to \textit{proactively} assign each operation of each order to specifically designated machines in anticipation of scrap and rework while minimizing the overall \textit{expected} tardiness. Accordingly, the overall layout of the \textit{proactive} scheduling consists of two attempts --- if scrap happens on a first attempt, the second attempt is initiated. These two attempts are captured probabilistically and are modeled within one problem formulation. The formulation developed by \cite{Liu21} is used as a basis, and stochastic elements due to scrap and rework are added, as explained in the following sections. Once a feasible schedule is obtained, the schedule corresponding to the first attempt is generated at the beginning of a shift: which job is to be assigned to what machine group and at what time. 
The integer-valued \underline{b}eginning and \underline{c}ompletion times for all operations within both attempts ((1) and (2)) are collectively denoted as vectors $b = \{\{b^1\},\{b^2\}\}$ and $c = \{\{c^1\},\{c^2\}\}$. 
After a schedule is obtained, a \textit{sequence of operations} to be processed by each machine group is inferred from the beginning times of the ``first attempt'' $\{b^1\}$ and the resulting sequence is passed on to the job-shop floor. After each shift, if there is no scrap, job-shop operations follow the aforementioned sequence. Otherwise, a rescheduling is triggered and a new sequence of operations is inferred. The short-solving times and high-quality schedules are thus especially important in the presence of uncertainties.  

\subsection{Objective Function}

To minimize the 
delays,
it is important to ensure the on-time shipment of processed orders. While it may not be possible to completely avoid delays, minimizing weighted tardiness can help mitigate their impact. Weighted tardiness is used as the objective function in the optimization problem, with the goal of minimizing the expected tardiness of all orders
\begin{align}
& o(c) = \label{eq1} \sum_{i=1}^I w_i \cdot \Bigg\{\prod_{j=1}^{J_i} (1-p^s_{(i,j)})\Bigg\} \cdot \max(c_{i,J_i}^1-d_i,0) + \\ \nonumber
& \sum_{i=1}^I w_i \! \cdot \!\sum_{j'=1}^{J_i} \Bigg\{\!\Bigg\{\!\prod_{j=1}^{j'-1} (1-p^s_{(i,j)}) - \prod_{j=1}^{j'} (1-p^s_{(i,j)})\Bigg\} \!\times\! \\
& (1-p^r_{(i,j')}) \!\cdot\! \max(c_{i,J_i,j',\bf{0}}^2-d_i,0)\!\Bigg\} + \nonumber \\ \nonumber
& \sum_{i=1}^I w_i \!\cdot\! \sum_{j'=1}^{J_i} \Bigg\{\!\Bigg\{\!\prod_{j=1}^{j'-1} (1-p^s_{(i,j)}) - \prod_{j=1}^{j'} (1-p^s_{(i,j)})\Bigg\} \!\times\! \\
& (p^r_{(i,j')}) \!\cdot\! \max(c_{i,J_i,j',\bf{1}}^2-d_i,0) \!\Bigg\}.  \nonumber 
\end{align}
The linearization of max operators follows the standard {\em special ordered set} procedure \cite{yan18, yan21}. 
Despite the complex appearance of the objective function within \eqref{eq1}, the important salient feature of it is \textit{additivity} in terms of jobs $i: o(c)=\sum_{i=1}^I o_i(c_i),$  which is efficiently exploited in Section \ref{sec:SolMethod}.

\subsection{Constraints}

To ensure the feasibility of the schedule, processing-time-requirement constraints (how long an operation needs to be processed), operation precedence constraints (to keep the order of operations as specified within a job), and machine capacity constraints (the expected number of operations to be processed on machines of a certain type cannot exceed the number of machines) are imposed at a modeling stage. 

\noindent \textbf{Beginning/Completion Time Constraints:} During Attempt $1,$ completion time $c_{i,j}^1$ of operation $(i,j)$ is expressed as:\footnote{Due to the discrete nature of time, one unit is subtracted. For instance, if an operation starts at 4 (e.g., 4:00 pm) with a processing time of 2, it finishes at 5 (e.g., 5:59:59 pm). If the next operation starts after the previous one's completion, it begins at time 6, leaving 2 units of time between the start times of subsequent operations.}   
\begin{equation}
 c_{i,j}^1 = \sum_{m \in E_{i,j}} \sum_{t = 1}^T (t+p_{i,j,m}) \cdot x_{i,j,m,t}^1 - 1, \forall(i,j),  \label{eq2}
\end{equation}
where $x_{i,j,m,t}^1$ is an indicator: $x_{i,j,m,\hat{t}}^1 = 1,$ if an operation $(i,j)$ begins at time $\hat{t}$~\footnote{Here $\hat{t}$ is used to distinguish a specific beginning time from $t,$ which is a dummy variable within summations.} on a machine group $m,$ and $0,$ otherwise; the above indicator variable can take the value of 1 for only one combination of $\{i,j,m,t\},$ therefore, the following constraint is imposed to satisfy Assumptions \ref{assume1}.\ref{a1} and \ref{assume1}.\ref{a3}: 
\begin{equation}
\sum_{m \in E_{i,j}} \sum_{t = 1}^T x_{i,j,m,t}^1 = 1, \forall(i,j). \label{eq3}
\end{equation}
Therefore, all the terms within \eqref{eq2} with the exception of the $\hat{t}^{\rm th}$
term are zero, and the part is completed at a time $\hat{t}+p_{i,j,m}.$ Accordingly, during Attempt 1, the beginning time $b_{i,j}^1$ of operation $(i,j)$ is expressed as: 
\begin{equation}
 b_{i,j}^1 = \sum_{m \in E_{i,j}} \sum_{t = 1}^T t \cdot x_{i,j,m,t}^1, \forall(i,j). \label{eq4}
\end{equation}

The completion time within the second attempt is modeled in a similar way: 
\begin{equation}
 c_{i,j,j',r}^2 \!= \!\!\!\!\sum_{m \in E_{i,j}} \!\!\sum_{t = 1}^T (t+p_{i,j,m}) \!\cdot\! x_{i,j,j',m,t,r}^2\! -\! 1, \!\forall(i,j,j',r). \label{eq5}
\end{equation}
The additional subscripts $j'$ and $r$ indicate that the completion time is computed for scenarios whereby a part is scrapped after operation $j'$ within the first attempt, and then the part is either reworked ($r=1$) or completely discarded ($r=0$). Analogously to \eqref{eq3} and \eqref{eq4}, the beginning times of the second attempt are formulated as follows:  
\begin{align}
&\sum_{m \in E_{i,j}} \sum_{t = 1}^T x_{i,j,j',m,t,r}^2 = 1, \forall(i,j,j',r),   \label{eq6}\\
&b_{i,j,j',r}^2 = \sum_{m \in E_{i,j}} \sum_{t = 1}^T t \cdot x_{i,j,j',m,t,r}^2, \forall(i,j,j',r).   \label{eq7}
\end{align}


\noindent \textbf{Operation Precedence:} Within each attempt, operation precedence constraints ensure that subsequent operation $j+1$ can only begin after the previous operation $j$ is completed:
\begin{equation}
 b_{i,j+1}^1 \geq c_{i,j}^1 + 1, \;\; b_{i,j+1,j',r}^2 \geq c_{i,j,j',r}^2 + 1, \forall (i,j,j',r).   \label{eq8}
\end{equation}


\noindent \textbf{Scrap/Rework Constraints:} Scrap/rework constraints ensure that the second attempt begins at Operation $1$ but only after operation $(i,j)$ leads to scrap ($r=0$) of part $i$: 
\begin{equation}
b_{i,1,j',0}^2 \geq c_{i,j'}^1 + 1, \forall(i,j').   \label{eq9}
\end{equation}
If a part can be reworked ($r=1$), then the operation $j'$ (second subscript) is repeated:  
\begin{equation}
 b_{i,j',j',1}^2 \geq c_{i,j'}^1 + 1, \forall(i,j').   \label{eq10}
\end{equation}

\noindent \textbf{Beginning of Shift Constraints:} It is a common practice to reschedule operations after the end of a shift to generate the schedule for the upcoming shifts. After part $i$ is scrapped after operation $j',$ the second attempt needs to be initiated at or after the beginning of future shifts. This is captured through the following constraints:
\begin{align}
 S \cdot \Big\lceil{\frac{c_{i,j'}^1}{S}\Big\rceil} + 1 &\leq  b_{i,j',j',0}^2,  S \cdot \Big\lceil{\frac{c_{i,j'}^1}{S} \Big\rceil} + 1 \leq  b_{i,1,j',1}^2, \forall(i,j'),   \label{eq10a}
\end{align}
where $S$ is the length of a shift. The ceiling operator is linearized after introducing integer variables $y_{i,j'}$ as: 
\begin{equation}
 \frac{c_{i,j'}^1}{S} \leq y_{i,j'} \leq \frac{c_{i,j'}^1}{S} + 1 - \varepsilon,   \label{eq10b}
\end{equation}
where $\varepsilon$ is a small positive number. 

\noindent \textbf{Expected Machine Capacity Constraints:} Machine capacity constraints ensure that at any point in time, the expected number of operations processed does not exceed the total number of machines within a group eligible to perform the corresponding operation: 
\begin{align}
& \sum_{(i,j) \in O_{m}} \prod_{j'=1}^{j-1} (1 - p^s_{{(i,j')}}) \cdot \left\{\sum_{t'=t-p_{i,j,m,1}: t'\geq 1}^t  x_{i,j,m,t'}^1\right\} + \label{eq11} \\
& \sum_{j', (i,j) \in O_{m}} \!\!\!\!\Bigg(\!\prod_{j=1}^{j'-1} (1\!-\!p^s_{{(i,j')}}) - \prod_{j=1}^{j'} (1\! -\! p^s_{(i,j')})\!\!\Bigg)  \!\times\! \nonumber\\
& (1- p^r_{{(i,j')}}) \!\cdot\! \left\{\sum_{t'=t-p_{i,j,m,1}: t' \geq 1}^t  \!\!\!\!\!x_{i,j,j',m,t',0}^2\!\right\} + \nonumber
\\
& \sum_{j', (i,j) \in O_{m}} \!\!\!\!\Bigg(\prod_{j=1}^{j'-1} \nonumber (1\!-\!p^s_{{(i,j')}}) - \prod_{j=1}^{j'} (1\! -\! p^s_{(i,j')})\!\!\Bigg) \!\times\! \nonumber\\
& p^r_{{(i,j')}}\! \cdot \!\left\{\sum_{t'=t-p_{i,j,m,1}: t' \geq 1}^t  \!\!\!\!\!x_{i,j,j',m,t',1}^2\!\right\}\!\leq\! C_{m}, \forall (m,t).  \nonumber 
\end{align}
Equation \eqref{eq11} ensures that at any time $t$ and for any machine group $m$, the expected number of parts processed does not exceed the number of machines $C_m$, and the expectation within \eqref{eq11} follows the similar logic as within the objective function \eqref{eq1}. The exception is the first term, thereby the upper limit of the summation is $j-1.$ For example, the probability of a part surviving \textit{up to} the last operation is $\prod_{j'=1}^{J_i-1} (1 - p^s_{{(i,j')}})$, while the probability of a part surviving \textit{after} the last operation is complete is $\prod_{j'=1}^{J_i} (1 - p^s_{{(i,j')}}),$ which is appropriate for calculating the tardiness as in \eqref{eq1}. 

For further convenience and compactness of notation, the objective function is denoted as $o(c) = \sum_{i=1}^I o_i(c_i)$ due to the additivity of \eqref{eq1} and the optimization problem becomes: 
\begin{equation}
\min_{b,c,x}\left\{\sum_{i=1}^I o_i(c_i), \;\; {\rm s.t.} \;\; \eqref{eq2}-\eqref{eq11}\right\}.  \label{eq12}
\end{equation}
The constraints \eqref{eq11} are additive in terms of parts $i$ and are converted to equality constraints through the use of nonnegative slack decision variables $z_{m,t}$, and are expressed as:  
\begin{equation}
\sum_{i=1}^I g_i(x_i) + z_{m,t} - C_{m} = 0.  \label{eq13}
\end{equation}
Here, $\sum_{i=1}^I g_i(x_i)$ compactly denotes the left-hand side of \eqref{eq11} with $x_i = \{\{x_{i,j,m,t}^1\},\{x_{i,j,j',m,t,r}^2\}\}$ denoting a vector containing all decisions with respect to operations and attempts. 
While all the variables within \eqref{eq1}-\eqref{eq11} are integer, because of the probabilities involved within \eqref{eq11}, the slack variables are modeled as continuous variables. In the following section, the solution methodology is presented.  

\section{Solution Methodology}
\label{sec:SolMethod}
Lagrangian Relaxation relies on the optimization of the dual function, which, in a general form can be written as: 
\begin{align}
 \max_{\lambda} \{q_\rho(\lambda): \lambda \in \mathbb{R}^M \times \mathbb{R}^T\},  \label{eq14}
\end{align}
with
\begin{align}
 & q_\rho(\lambda) = \label{eq15} \\
 & \min_{b,c,x,z} \Bigg\{L_\rho(c,x,z,\lambda), (b_i,c_i,x_i) \in \mathcal{F}_i, i = 1,\dots,I, z \geq 0. \Bigg\},  \nonumber
\end{align}
where 
\begin{align*}
L_\rho(c,x,z,\lambda) \equiv &\sum_{i=1}^I o_i(c_i) \!+\! \lambda^{\rm T} \!\cdot\! \left(\sum_{i=1}^I g_i(x_i) + z_{m,t} - C_{m}\right)\\
& \!+\! \rho \!\cdot \!\left\|\sum_{i=1}^I g_i(x_i) + z_{m,t} - C_{m}\right\|_1\!,    
\end{align*}
is the ``absolute-value'' Lagrangian function \cite{Bragin19}.
The feasible set $\mathcal{F}_i$ for each job $i$ is delineated by constraints \eqref{eq2}-\eqref{eq10}. 
Lagrangian multipliers $\lambda$ (``dual'' variables) are the decision variables with respect to the dual problem \eqref{eq14} and it is assumed that the set of optimal solutions $\Omega = \{\lambda \in \mathbb{R}^M \times \mathbb{R}^T \;\; | \;\; q(\lambda) = q(\lambda^*)\}$ is not empty. The minimization within \eqref{eq15} with respect to $\{b,c,x,z\}$ is referred to as the ``relaxed problem''. Due to integer variables $\{b,c,x,z\}$, the dual function \eqref{eq14} is non-smooth with facets (each corresponding to a particular solution to the relaxed problem within \eqref{eq15}) intersecting at ridges whereby derivatives of $q(\lambda)$ exhibit discontinuities; in particular, the dual function is not differentiable at $\lambda^{*}$. 

To maximize the dual function, ``Surrogate Level-Based Lagrangian Relaxation'' \cite{bragin2022surrogate} is chosen. The dual function \eqref{eq14} is maximized by updating Lagrange multipliers $\lambda$ by taking a series of steps $s^k$ along ``surrogate'' sub-gradient directions $\sum_{i=1}^I g_i(\tilde{x}_i^k) + \tilde{z}_{m,t}^k - C_{m}$ (levels of constraint violation) as 
\begin{equation}
\lambda^{k+1} = \left[\lambda^k + s^k \cdot \left(\sum_{i=1}^I g_i(\tilde{x}_i^k) + \tilde{z}_{m,t}^k - C_{m}\right)\right]^+,   \label{eq16}
\end{equation}
where $[\cdot]^+$ is a projection operator onto a positive orthant $\left\{\lambda | \lambda \geq 0\right\}$. Following \eqref{SLBLR}, the step-sizes are set as: 
\begin{equation}
 s^k = \!\zeta\! \cdot\! \gamma\! \cdot\! \frac{\overline{q}_{k} - L(\tilde{b}^k,\lambda^k)}{\left\|\sum_{i=1}^I g_i(\tilde{b}_i^k)\! +\! \tilde{s}_{m,t}^k \!-\! M_{m}\right\|^2_2}, \gamma < 1, \zeta < 1.  \label{eq17}
\end{equation}
Here the ``$\sim$'' is used to distinguish solutions $\{b^k,x^k,z^k\}$ that are obtained by solving optimally the relaxed problem (minimization within \eqref{eq15}) from subproblem solution $\{\tilde{b}^k,\tilde{x}^k,\tilde{z}^k\}$. Subproblems are formulated by forming a group of parts $\mathcal{I}$ and optimizing the relaxed problem with respect to the associated variables while keeping decision variables not belonging to a subset $\mathcal{I}$ fixed as:  

\begin{align}
\min_{b_i,c_i,x_i,z; i \in \mathcal{I}} &\left\{\sum_{i \in \mathcal{I}} o_i(c_i) + \lambda^{\rm T} \cdot \left(\sum_{i \in \mathcal{I}} g_i(x_i) + z_{m,t}\right)+ \right.\nonumber \\
& \!\!\!\!\!\!\!\!\!\!\!\!\!\!\!\!\!\!\!\!\!\!\!\!\!\!\!\!\!\!\!\!\!\!+\left.\rho \cdot \left\|\sum_{i \in \mathcal{I}} g_i(x_i) + \sum_{i \notin \mathcal{I}} g_i(x_i^{k-1}) + z_{m,t} - C_{m}\right\|_1\right\}.  \label{eq20a}
\end{align}
\noindent The above minimization involves piecewise linear penalties ($l_1$ norms), that efficiently penalize constraint violations and are exactly linearizable through the use of \textit{special ordered sets}, thereby enabling the use of MILP solvers. 

\textbf{Heuristics.} Heuristics are generally the last, yet, necessary step required to obtain feasible solutions. Feasible solutions are obtained by repairing subproblem solutions and the feasible solution search is initiated after the norm of constraint violations is “sufficiently small.” For the JSS problem, the following simple rule along the lines of a ``local search'' is used: \textit{adjust subproblem beginning times of operations $(i,j)$ by no more than $\delta_j$ time periods}. This rule is operationalized by solving the original problem (to enforce feasibility) subject to extra restrictions to model the elements of the ``local search'' as follows:
\begin{equation}
\min_{b,c,x,z} \;\;\eqref{eq1}, \;\; {\rm s.t.}\;\; \eqref{eq2}-\eqref{eq11}, \;\;- \delta_j \leq b_{i,j} - b_{i,j}^k \leq \delta_j. 
\label{eq21}
\end{equation}
The above problem is then solved by a MILP solver. 
The steps of the SLBLR method are summarized in the algorithm below. 

\begin{algorithm}{Surrogate ``Level-Based'' Lagrangian Relaxation}\label{alg:SLBLR}

\begin{itemize}
\item[] \textbf{Step 1: Initialization.} Initialize $\lambda^0,$ $s^0,$ $c^0.$
\item[] \textbf{Step 2: Subproblem Solving.} Solve a subproblem \eqref{eq20a}. 
\item[] \textbf{Step 3: Stepsize Update.}  Use \eqref{eq17} to update stepsizes.
\item[] \textbf{Step 4: Multiplier Update.}  Use \eqref{eq16} to update multipliers.
\item[] \textbf{Step 5: Penalty Coefficient Update.} Update $\rho^{k+1} = \rho^{k} + \beta, \beta > 0.$\footnote{When solving continuous problems, penalty coefficients are typically updated by using a \textit{multiplicative} constant $\beta$ as $\rho^{k+1} = \rho^{k} \cdot \beta, \beta > 1$, for example, within the Method of Multipliers \cite{bertsekas97}. However, an $additive$ constant $\beta$ is adopted here.}
\item[] \textbf{Step 6: Constraint Violation Check.} Check the Levels of Constraint Violation. If $\left\|\sum_{i=1}^I g_i(x_i) + z_{m,t} - C_{m}\right\| \leq \epsilon$, go to Step 7, otherwise, go to Step 2. 
\item[] \textbf{Step 7: Feasible Solution Search.} Solve \eqref{eq21}. If the CPU time limit is not reached, go to 2, otherwise, terminate with the best feasible solution obtained. 
\end{itemize}
\end{algorithm}

The efficiency of the overall approach is discussed with respect to numerical
case studies in the following section. 
\section{Numerical Testing Results}\label{sec:numerical_examples}

The solution methodology was implemented in an open-source Julia package, Jobshop.jl available at \url{https://github.com/PSORLab/Jobshop.jl} \cite{JobshopPackage}. Each subproblem in the numerical examples was solved using a 64-bit CPLEX 12.10.0 optimizer and tested using a CPU (48 threads) of Intel(R) Xeon(R) E-2286M CPU @ 2.99 Hz, 192 GB of RAM, and Windows 10. Example 1 is to demonstrate the computational efficiency of the method, it is an instance with 20 jobs, 5 operations per job, and 5 machines, each specializing in a specific operation. Example 2 is to demonstrate the scalability of the SLBLR method with respect to the increase in the number of jobs as well as in duration of the processing times.

\subsection{Example 1. Instance with 20 Jobs and Short Processing Times.} Within this example, a base-case instance with 20 jobs, 5 operations per job, and 5 machines, each designated to perform one specific operation is considered. Processing times are generated randomly based on discrete uniform distribution $U[1,5]$. The data used in this example are shown in Table \ref{Table:example1} in the Appendix. 

\begin{figure}[ht]
		\centering
 \includegraphics[trim = 0 0 0 15, scale=0.55, angle=0]{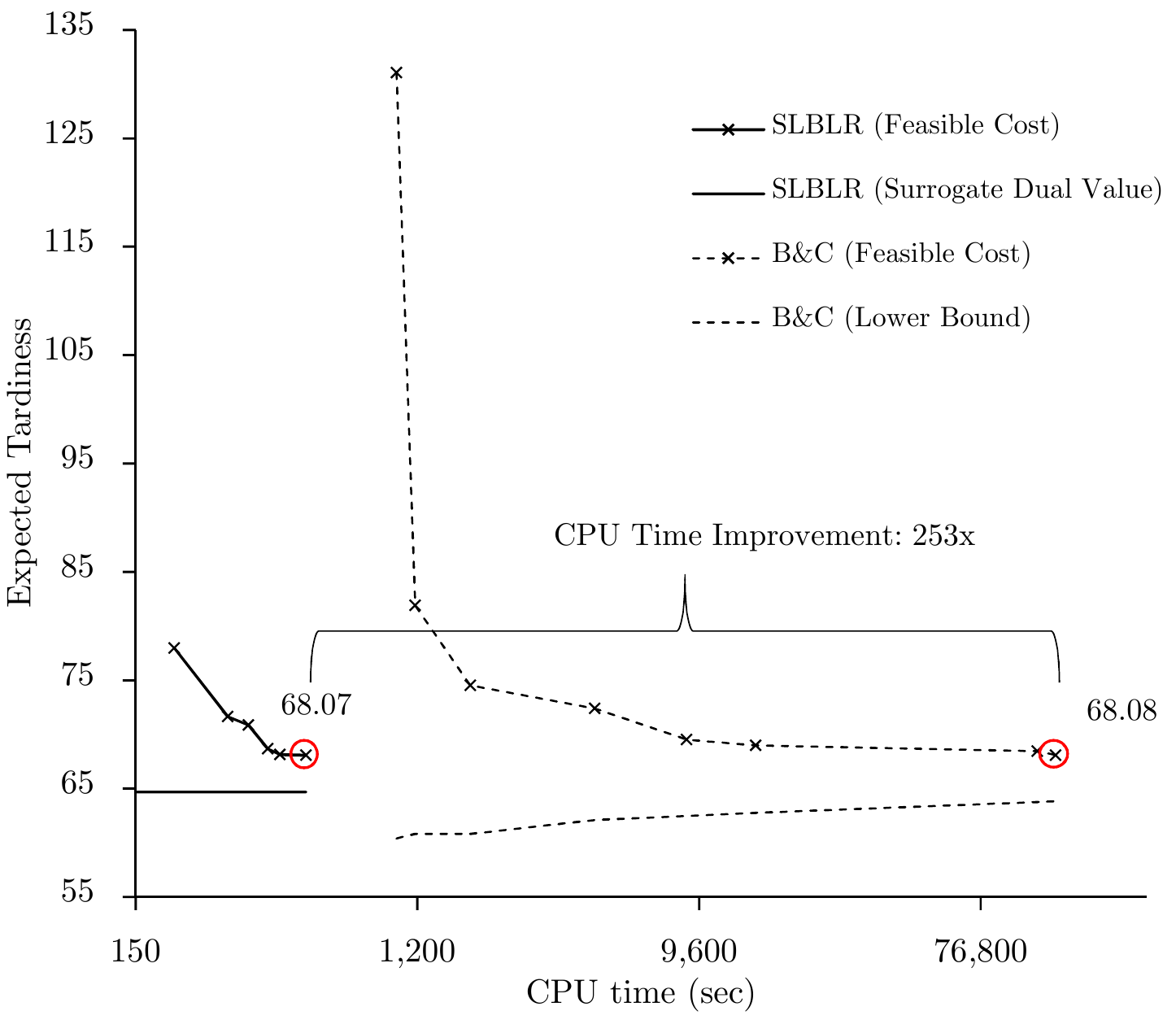}	\caption{The results for Example 1 are illustrated. The SLBLR method
 performs approximately 250 times faster than the B\&C algorithm.}
 \label{fig:ex1_results}
\end{figure}

The due dates are $d = \{$15, 25, 32, 36, 21, 27,26,	13,	29,	12,	35,	31,	19,	24,	33, 23,	18,	21,	17,	22$\}.$ The machine capacities are $M = \{2, 3, 2, 2, 3\}.$ The scrap rate is 5\% and the rework rate is 20\%. The priorities are all chosen to be equal $w_i = 1, i = 1,\dots, I$. The results of the SLBLR method and standard B\&C are shown in Figure \ref{fig:ex1_results}.

As illustrated in Figure \ref{fig:ex1_results}, the SLBLR method is more than two orders of magnitude faster than standard B\&C.

To test robustness, the base case is modified by re-generating processing times based on the discrete uniform distribution $U[1,5]$ as well as by generating due dates based on discrete uniform distribution $U[10,40].$ Accordingly, 5 new instances are generated. Because of the complexity of the problem, the stopping criterion is 10\% of the duality gap. The results are demonstrated in Table \ref{Ex3}.  
\begin{table*}[!ht]
\small
\caption{Robustness results for Example 1.}
\label{Ex3}
\centering
\scalebox{1}{
\begin{tabular}{|c|c c c c|c c c|c|}
\hline
 &  \multicolumn{4}{c|}{B\&C} & \multicolumn{3}{c|}{SLBLR} &\multicolumn{1}{c|}{}\\
\hline
Case	&	Feas.  	&	Lower 	&	Gap (\%)	&	CPU 	&	Feas. 	&	Gap (\%)	&	CPU &	Improvement	\\
	&	Cost	&	Bound 	&		&	time (sec)	&	Cost	&		&	time (sec)	&		\\
\hline
1	&	133.50	&	120.27	&	9.91\%	&	12189.88	&	131.45	&	8.51\%	&	401.20	&	30.38	\\
2	&	158.89	&	143.25	&	9.84\%	&	78690.8	&	158.82	&	9.80\%	&	313.17	&	251.26	\\
3	&	137.05	&	123.34	&	10.00\%	&	78788.47	&	135.31	&	8.85\%	&	843.47	&	93.40	\\
4	&	150.31	&	135.31	&	9.98\%	&	6610.47	&	149.83	&	9.69\%	&	322.04	&	20.52	\\
5	&	118.52	&	106.98	&	9.74\%	&	41301.77	&	118.28	&	9.56\%	&	404.07	&	102.21	\\
\hline
	&		&		&	Avg. time:	&	43516.27	&		&	Avg. time:	&	456.79	&	\textbf{95.26}	\\
	&		&		&	St. Dev.:	&	34747.87	&		&	St. Dev.:	&	220.32	&		\\

\hline
\end{tabular}
}

\end{table*}
 
Table \ref{Ex3} demonstrates that the CPU time improvement ranges from 20.52 times to 251.26 times, for an average of 95.26 times, which is almost two orders of magnitude improvement.

\subsection{Example 2. Instances with 20 and 100 Jobs and Long Processing Times.} 
\begin{figure}[ht]
		\centering
 \includegraphics[trim = 0 0 0 0, scale=0.55, angle=0]{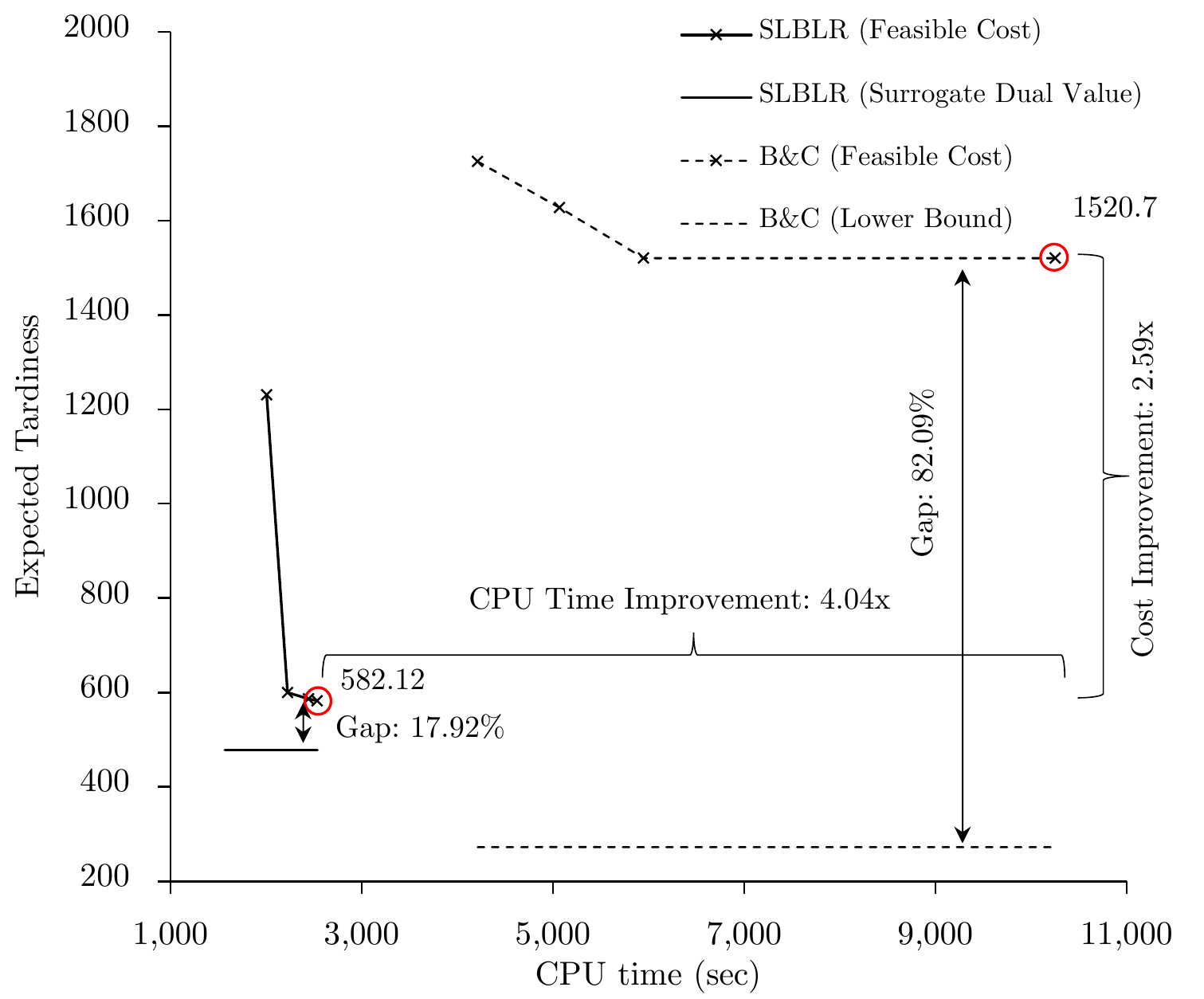}	\caption{The results for Example 2 are illustrated. The SLBLR method
 performs approximately 4 times faster than the B\&C algorithm and 2.59 times more accurately.}
 \label{fig:ex2_results}
\end{figure}

\noindent \textbf{Results for a 20-Job Case.} The base case study from Example 1 is modified by increasing due dates by a factor of 10 and by generating processing times by using uniform discrete distribution $U[1,50].$ The results are shown in Figure \ref{fig:ex2_results}. 

Because of the associated increase of the time periods required to accommodate much-increased processing times, the CPLEX is out of memory after 10,255 seconds. In contrast, the SLBLR method is 4 times faster and 2.59 times more accurate. 

\noindent \textbf{Results for a 100-Job Case.} The instance is created by generating processing times based on the discrete uniform distribution $U[1,50]$ as well as by generating due dates based on discrete uniform distribution $U[100,400].$ For this case, no feasible solution was obtained by CPLEX after 15,129.42 sec. In contrast, a feasible solution with a corresponding cost of 9697.59 is found by SLBLR after 3618.53 seconds. 



\section{Conclusions, Limitations, and Future Work.}
\label{sec:conclusion}
This study presented a proactive scheduling approach that minimizes the expected tardiness in the presence of scrap and rework in job-shop manufacturing environments. The major contributions and takeaways of this study can be summarized as follows:

\begin{enumerate}
\item \textbf{Formulation.} The use of Markov processes to model probabilistic defects in manufacturing production, and the formulation of an optimization problem that minimizes expected tardiness.
\item \textbf{Online Capabilities of SLBLR.} Drawing on the geometric potential inherited from the Polyak formula, the SLBLR method can effectively operate online and accommodate unexpected or urgent orders by continuously updating multipliers. This capability enhances the resilience of manufacturing operations.
\item \textbf{Consistent Improvement Capabilities.} The SLBLR method's ability to systematically improve solutions through the update of multipliers.
\item \textbf{Optimization Enhancement through QC and ML Integration.} The SLBLR method's adaptability to different optimization techniques facilitates seamless integration with quantum computing (QC) and machine learning (ML). This combination holds the potential to substantially improve subproblem-solving and feasible solution quality. Recent advancements have shown promising results in tackling job-shop scheduling challenges within the automation community, highlighting the value of merging Lagrangian Relaxation and Machine Learning methods \cite{liu2023integrating}.
\item \textbf{Broader Impact.} The potential for significant reductions in tardiness and improved downstream supply chain operations management.

\end{enumerate}

While the results and future directions are promising, some limitations should be acknowledged. First, the computational effort can be considerable for large-scale problems, as seen in the 100-job case. Second, the model assumes all jobs are available at the beginning of each shift, which may not always be accurate. Lastly, the model does not account for sequence-dependent setup times, which can significantly affect manufacturing operations.

Future research can focus on overcoming these limitations by developing more efficient algorithms for large-scale problems, incorporating real-time job availability, and considering setup times in the optimization formulation. Additionally, exploring the integration of the SLBLR method with emerging technologies like artificial intelligence and Industry 4.0 could lead to more comprehensive and effective manufacturing operations management.







\section*{Disclaimer}
This paper was prepared as an account of work sponsored by an agency of the United States Government. Neither the United States Government nor any agency thereof, nor any of their employees, makes any warranty, express or implied, or assumes any legal liability or responsibility for the accuracy, completeness, or usefulness of any information, apparatus, product, or process disclosed or represents that its use would not infringe privately owned rights. Reference herein to any specific commercial product, process, or service by trade name, trademark, manufacturer, or otherwise does not necessarily constitute or imply its endorsement, recommendation, or favoring by the United States Government or any agency thereof. The views, opinions, and/or findings contained in this paper are those of the authors and should not be interpreted as representing the official views or policies, either expressed or implied, of the United States Government or any agency thereof, including the Air Force Research Laboratory, the United States Air Force, and the Department of Defense.

\section*{Funding}

This material is based upon work supported by the U.S. Department of Energy’s Office of Energy Efficiency and Renewable Energy (EERE) under the Advanced Manufacturing Office Award No. DE-EE0007613.
We also gratefully acknowledge the Air Force Research Laboratory, Materials and Manufacturing Directorate (AFRL/RXMS) for support via Contract  No. FA8650-20-C-5206. This work is also supported in part by the US NSF under award ECCS-1810108.

\bibliographystyle{IEEEtran}
\bibliography{bibliography.bib}

\appendix
\section{Example Data}
\setcounter{table}{0}
\begin{table}[!ht]
\caption{Processing time $p_{i,j}$ data for job $i$ and operation $j$ for Example 1 are contained in this table.}

\label{Table:example1}
\centering
		{\def\arraystretch{0.2}  
\begin{tabular*}{0.34\textwidth}{@{\hspace{0pt}}ccc|ccc|ccc}
\hline
\hline
$i$ & $j$ & $p_{i,j}$ & $i$ & $j$ & $p_{i,j}$ & $i$ & $j$ & $p_{i,j}$  \\
\hline
1	&	1	&	1	&	2	&	1	&	2	&	3	&	1	&	3	\\
1	&	2	&	2	&	2	&	2	&	1	&	3	&	2	&	2	\\
1	&	3	&	2	&	2	&	3	&	2	&	3	&	3	&	3	\\
1	&	4	&	1	&	2	&	4	&	2	&	3	&	4	&	2	\\
1	&	5	&	1	&	2	&	5	&	4	&	3	&	5	&	2	\\
4	&	1	&	1	&	5	&	1	&	3	&	6	&	1	&	1	\\
4	&	2	&	3	&	5	&	2	&	2	&	6	&	2	&	3	\\
4	&	3	&	2	&	5	&	3	&	5	&	6	&	3	&	2	\\
4	&	4	&	2	&	5	&	4	&	1	&	6	&	4	&	3	\\
4	&	5	&	1	&	5	&	5	&	4	&	6	&	5	&	2	\\
7	&	1	&	2	&	8	&	1	&	1	&	9	&	1	&	1	\\
7	&	2	&	5	&	8	&	2	&	5	&	9	&	1	&	2	\\
7	&	3	&	4	&	8	&	3	&	1	&	9	&	1	&	2	\\
7	&	4	&	1	&	8	&	4	&	5	&	9	&	1	&	5	\\
7	&	5	&	1	&	8	&	5	&	1	&	9	&	1	&	3	\\
10	&	1	&	1	&	11	&	1	&	5	&	12	&	1	&	3	\\
10	&	1	&	3	&	11	&	2	&	2	&	12	&	2	&	1	\\
10	&	1	&	2	&	11	&	3	&	3	&	12	&	3	&	3	\\
10	&	1	&	3	&	11	&	4	&	4	&	12	&	4	&	1	\\
10	&	1	&	1	&	11	&	5	&	5	&	12	&	5	&	1	\\
13	&	1	&	5	&	14	&	1	&	2	&	15	&	1	&	3	\\
13	&	2	&	1	&	14	&	2	&	4	&	15	&	2	&	1	\\
13	&	3	&	3	&	14	&	3	&	2	&	15	&	3	&	1	\\
13	&	4	&	1	&	14	&	4	&	1	&	15	&	4	&	3	\\
13	&	5	&	5	&	14	&	5	&	1	&	15	&	5	&	5	\\
16	&	1	&	3	&	17	&	1	&	3	&	18	&	1	&	1	\\
16	&	2	&	3	&	17	&	2	&	1	&	18	&	2	&	2	\\
16	&	3	&	2	&	17	&	3	&	2	&	18	&	3	&	5	\\
16	&	4	&	5	&	17	&	4	&	3	&	18	&	4	&	2	\\
16	&	5	&	3	&	17	&	5	&	1	&	18	&	5	&	3	\\
19	&	1	&	1	&	20	&	1	&	2	&		&		&		\\
19	&	1	&	4	&	20	&	1	&	4	&		&		&		\\
19	&	1	&	4	&	20	&	1	&	4	&		&		&		\\
19	&	1	&	2	&	20	&	1	&	2	&		&		&		\\
19	&	1	&	5	&	20	&	1	&	1	&		&		&		\\

\hline
\hline
\end{tabular*}
}
\end{table}

\end{document}